\documentclass[a4paper]{article}

\usepackage[english]{babel}
\usepackage[utf8]{inputenc}
\usepackage{amsmath}
\usepackage{graphicx}
\usepackage[colorinlistoftodos]{todonotes}
\usepackage{psfrag,epsf}
\usepackage{enumerate}
\usepackage[numbers]{natbib}
\usepackage{url} 
\usepackage{amsthm,mathrsfs,amsfonts}	
\usepackage{enumitem}
\usepackage{bbm}
\newcommand*{\QEDB}{\hfill\ensuremath{\square}}%
\newcommand{\X}{ \textbf{X}}

\def\defeq{\mathrel{\mathop:}=}
\newcommand{\eqdef}{\mathrel{\mathop:}=}

\numberwithin{equation}{section}
\newtheorem{theorem}{Theorem}[section]

\newtheorem{remark}{Remark}

\newtheorem{corollary}{Corollary}

\def\defeq{\mathrel{\mathop:}=}

\newcommand{\ccs}{\mbox{$\stackrel{c.s}{\longrightarrow}$}}

\newcommand{\M}{\mbox{$\mathcal{M}$}}

\newcommand{\I}{\mathbbm{1}}
\newcommand{\R}{\mbox{$\mathbb{R}$}}

\usepackage{graphicx} 
\usepackage{subfig}
\usepackage{float}		
\usepackage{color}
\usepackage[normalem]{ulem}
\title{\bf Connecting pairwise geodesic spheres by depth: DCOPS}

\author{Ricardo Fraiman \thanks{rfraiman@cmat.edu.uy}\hspace{.2cm}\\
    Centro de Matem\'atica, Facultad de Ciencias, \\
     Universidad de la Rep\'ublica, Uruguay.\\
   Fabrice Gamboa \\
    Institut de Math\'ematiques de Toulouse, France.\\
 and \\ 
Leonardo Moreno \\ 
Departamento de M\'etodos Cuantitativos, FCEA, \\
Universidad de la Rep\'ublica, Uruguay.}

\date{2018-04-27}

\begin{document}
\maketitle

\begin{abstract}
We extend the classical notion of the spherical depth in $\mathbb{R}^k$, to the important setup of data on a  Riemannian manifold. We show that this notion of depth satisfies a set of desirable properties. For the empirical version of this depth function both uniform consistency and the asymptotic distribution are studied. Consistency is also shown for functional data. The behaviour of the depth is illustrated through several examples for Riemannian manifold data.
\end{abstract}

\noindent%

\section{Introduction}	

Depth concepts have become a fundamental tool in modern statistics.  Over the last three decades, data depth has emerged as a powerful concept, which leads to a center-outward ordering of sample points in multivariate data and, more recently, in functional data. Besides providing a robust notion of multivariate center (median) and trimmed means, applications to many other statistical problems have recently been developed. They have been successfully applied to supervised and unsupervised learning, hypothesis testing, and outliers detection; see, for instance, the seminal paper by \cite{liu1999}. The notion of depth was first introduced by \cite{tukey1975} for bivariate data sets and it was extended to higher dimensions by \cite{donoho1992}.  Several other notions were later introduced, such as \textit{convex hull peeling depth}, \cite{barnett1976},\textit{Oja depth} \cite{oja1983},  \textit{simplicial depth} \cite{liu1990}, \textit{spatial depth} \cite{vardi2000}, and spherical depth \cite{elmore2006} among others. In particular, several notions have been introduced in the functional data setting in the last few years.

Time has shown that half--space, simplicial, spatial depth and spherical depth are fundamental and their properties have been widely studied. They share a list of desirable properties  that were introduced in \cite{liu1990} for the simplicial depth and studied by \cite{serfling2000} for general depth notions. These desirable properties are
\begin{itemize}
\item [P1)]  \it Invariance with respect to a group of transformations. Affine invariance or orthogonal invariance. 

\rm P11: \it Affine invariance: \rm  The depth should not depend on the coordinate system and, in particular, on the scales of the underlying measurements.

P12 \it Orthogonal invariance: \rm The depth should not depend on the coordinate system and the global scale (the procedure is equivariance under location, unitary operators and homogeneous scale transformations). 
\item [P2)] \it Maximality at center: \rm For a distribution having a uniquely defined ``center", such as the point of symmetry with respect to some notion of symmetry (invariance under a suitable family of transformations like angular, central, elliptical or spherical symmetry), the depth should attain a maximum value at this center. Observe that all of these notions of symmetry coincide with the usual notion of symmetry for one dimensional data. 
\item [P3)] \it Monotonicity relative to deepest point. When a point moves away from the ``deepest point'', through  a ray,  the depth decreases. \rm
\item [P4)] \it Vanishing at infinity. \rm 
\end{itemize}

Serfling (see \cite{serfling2006}) mentioned two extra important desirable properties:

``\textit{Dimensionality}. Typically, data in $\mathbb{R}^k$ has structure of lesser dimension than
the nominal $k$.''

``\textit{Computional  complexity}. The feasibility of computation is a function of the
size n and and dimension d of the data, with practical limitations arising in the
case of higher $k$.''

Nevertheless, with the exception of spherical depth the effective calculation for data in high dimensional spaces (or even moderate dimension) of these classical depths is a hard computational problem. In particular this property is very important in the modern image data analysis, see \cite{pizer2017}.  Another important issue is that nowadays many statistical problems have to deal with data that are not in finite dimensional Euclidean spaces, such as functional data or data in some metric spaces, in particular data in $\mathbb{R}^k$ with structure of lesser dimension than
the nominal $k$.
The case of the statistical analysis of shape distributions with applications in medical diagnostic, morphometrics and for the analysis based on planar views of 3-D objects is of particular importance. The setup where data are on Riemannian manifolds is the appropiate one for the statistical analysis.

Our main objective here is to extend spherical depth  to more general spaces, particularly for data on Riemannian manifolds. We will refer in what follows to spherical depth as DCOPS. This  statistical frame arises in several important applications; see, for example,  \cite{patrangenaru2015}, \cite{bhattacharya2008} and \cite{pennec2006}. In particular, defining centrality measures on Riemannian manifolds is an important issue; see  \cite{fletcher2009} and \cite{barbaresco2013}.

In FDA different notions of depth functions have been proposed and, as mentioned in \cite{cuevas2014}, is far from being an exhausted topic.  We show that, under weak conditions, the empirical depth proposed in \cite{elmore2006} is also consistent for functional data spaces.

For both problems, data on Riemannian manifolds as well as functional data  new technical challenges appear. We have to deal with geodesic balls and its properties for Riemannian manifolds to derive the main properties as well as the asymptotic behaviour of DCOPS, while for FDA we need to make use of  some fundamental results of \cite{billingsley1967} that provide necessary and sufficient conditions for uniform convergence over a given family of functions. Further, the important case where data are random positive definite matrices is studied in  detail.

Two approaches have been considered to model  data  on a manifold. Either  embed the manifold in an ambient Euclidean space (extrinsic approach) or only use the intrinsic properties of the Riemannian manifold (intrinsic approach).  For instance, for  the regression problem, we could mention \cite{lin2016} and \cite{zhu2009}, respectively. In this paper we follow the second approach. The first advances in the intrinsic approach (also in regression) was given by \cite{pelletier2006}  using the exponential map locally.
The computational problem of simplicial depth derives from the fact that the difficulty of checking if a point belongs to a simplex increases considerably with the dimension of the space. One way to circumvent this problem is to consider spherical depth,  a depth that is based on pairwise connecting spheres, introduced in \cite{elmore2006}.  More precisely, for the Euclidean space, it was defined  as follows:
given two points $x,y \in \mathbb R^k$, we denote by $B_{xy}$ the closed ball whose diameter is the segment $\overline{xy}$; that is, the ball  centered at the midpoint between $x$ and $y$ (which is denoted by $\frac{xy}{2}$ in what follows) and with radius $d(x,y)/2$. Then, we define the population depth by connect pairwise spheres (DCOPS) of a random vector $X$ as
\begin{equation}
\begin{aligned}
BD(x) \eqdef P\left( x \in B_{X_1 X_2} \right),
\end{aligned}
\label{defprofundidad_euc}
\end{equation}
where $x \in \mathbb R^k$, and  $X_1$ and $X_2$ are two independent random vectors with the same distribution as  $X$.
Given an iid random sample  $ \{X_1, \ldots, X_n\}$ with the same distribution as $X$,  the empirical version of $BD$ is given by the following $U$-statistics of order $2$
\begin{equation}
\widehat{BD}_n(x) \defeq \frac{1}{\binom{n}{2}}\sum_{1 \leq i_1 < i_2 \leq n} \I_{B_{X_{i_1} X_{i_2}}}(x)\defeq U^n_2 \left( D_x\right),
\label{estimadorprofundidad_euc}
\end{equation}
where
\begin{equation}
\label{D}
D_x \defeq \{ (x_1,x_2) \in \mathbb{R}^k \times \mathbb{R}^k: B_{x_1 x_2} \ni x \} \textrm{ and we denote } \mathcal{D} \defeq \left\{ D_x \right \}_{x \in \mathbb{R}^k}.
\end{equation}

Notice that $BD$ coincides with simplicial depth when $k=1$. Also, it is well defined for data in any separable Hilbert space $\mathcal H$, where it also has an equivalent  geometric interpretation given by
\begin{equation}\label{defpob}
BD(x)= P \left( \langle X_1-x, X_2-x \rangle \leq 0 \right),  \, \left(x \in \mathcal{H}\right)
\end{equation}
where $X_1$ and $X_2$ are two independent random elements in $\mathcal H$ with the same distribution as $X$, and we denote by  $\langle \cdot, \cdot \rangle$ the inner product on $\mathcal H$.
 The computational complexity of this notion is just of order $O(k \times n^2)$, linear on the dimension $k$ and quadratic on the sample size $n$. Some authors claim that a shortcoming of spherical depth is that it is only orthogonal invariant and not affine invariant. However, in our setup (in Riemannian manifolds), it only makes sense to be orthogonal invariant. 
 \\ 
This paper is organized as follows. In Section \ref{setup} we introduce the setup for random element taking values on a Riemannian manifold. In Section \ref{dcops} we prove the basic properties of DCOPS for the case where the data take values on a Riemannian manifold. In Subsection \ref{pdcops}  we  obtain results on the uniform convergence and on the asymptotic distribution of the empirical version of the depth function. The case of FDA is studied in Subsection \ref{fda}.  Section \ref{simulated} provides some simulated examples for  non--Euclidean spaces. In particular, we consider the important case of the space of positive definite  random matrices. Some technical proofs are postponed to the Appendix.

\section{Riemannian manifolds spaces}
\label{setup}
In this section we consider the case where the random elements take values on a Riemannian manifold. 
\subsection{The setup}
To begin with, we give some tools and facts on Riemannian geometry, that are useful in our manuscript. A \textit{Riemannian metric} $g$ on a manifold $\M$ defines for every point $p \in \M$  the scalar product (denoted by $g_p(\cdot,\cdot)$). It smoothly depend on the point $p$,  of tangent vectors in the tangent space $T_p \M$. The \textit{Riemannian manifold} $(\M,g)$ is a manifold that is equipped with the Riemannian metric  $g$. 
Given $p \in \M$, for every $v \in T_p \M$ the Riemannian norm of $v$ is given by $\Vert v \Vert \defeq \sqrt{g_p(v,v)}$.
 If  $\gamma: [a, b] \rightarrow \M$ is a continuously differentiable curve in the Riemannian connected manifold $\M$ and  $L(\gamma)$ is the length of the curve $\gamma$,  then the \textit{induced distance}  $d(x, y)$ between the points $x$ and $y$ of $\M$ is defined as
  \begin{align*}
 d(x,y)= \\
 &\hspace{-1.7 cm}= \inf \{ L(\gamma) : \gamma \, \textrm{is a continuously differentiable curve joining $x$ and $y$}\}.
  \end{align*}

The letter $d$ will stand for the Riemannian distance on $\M$  with respect to the $g$ metric. A \textit{geodesic} (with speed $s \in \mathbb{R}^{+}_0$) is a smooth map $\alpha:I \rightarrow \M$, where $I$ is an interval, such that $\Vert \alpha^{'}(t) \Vert= s$  for all $t \in I$ and which is ``locally length minimizing'' . For any $p \in \M$ and a vector $v$ in the tangent space $T_p \M$, there exists a unique geodesic  $\alpha_{(p,v)}(t)$ starting from that point with this tangent vector $v$. The \textit{exponential map} is the map $exp_p$ given by $\exp_p (v) = \alpha_{(p,v)}(1)$. 
\begin{remark}
The geodesic that joins two points need not be globally length--minimizing and it is not necessarily unique.
\end{remark}
 The \textit{cut locus} of $p$  in $\M$  is the maximal domain where the exponential map is a diffeomorphism. That is, the cut locus of $p$ in the tangent space is defined to be the set of all vectors $v \in T_{p} \M$ such that  $\exp _{p}(tv)$ is a minimizing geodesic for  $0 \leq t \leq 1$  but fails to be one for $0 \leq t \leq 1+ \epsilon$, $\epsilon >0$ . The \textit{cut locus of $p$ in $\M$}, denoted  by $C_{\M}(p)$, is defined as the image of the cut locus of  $p$ in the tangent space under the exponential map at $p$. The \textit{injectivity radius of $p$} is the maximal radius of centered balls on which the exponential map is a diffeomorphism. The \textit{injectivity radius of the manifold}  $r_{iny}$ is the infimum of the injectivity over the manifold.

Let $(\M,g)$ be a connected and orientable Riemannian manifold, see \cite{docarmo1992}, page $18$. 
We will assume that  $(\M,d)$ is a complete  separable metric space. Since $(\M,d)$ is complete, the Hopf--Rinow theorem 
(see \cite{docarmo1992}, p.146) implies that for any pair of points $p,q \in \M$  there exist at least one geodesic path in $\M$ connecting $p$ and $q$. 
   
  If the manifold $\M$ fulfills the assumptions of the Cartan--Hadamard Theorem (see \cite{petersen2006}, p.162), that is, if $\M$ is a simply connected, complete Riemannian manifold with non--positive curvature (a Hadamard manifold), then the geodesic is unique.
  
  Let $X$ be a random element taking values in $\M$, with distribution $P$. If some of these conditions for the uniqueness of the geodesic are not fulfilled, then we will need to assume some more condition on $P$ that implies that there is a unique geodesic between each pair of points $p,q \in \M$ ($P \times P$ almost surely). More precisely, we will assume that the random element $X$ has a density $f_X$ with respect to the volume measure $d \nu(x)$ on $\M$, (which exists since $\M$ is orientable, see  the last section of \cite{folland2013}) fulfilling the following condition:
Given $q \in \M$ let 
\begin{equation}
A_q\defeq \left\{x \in \M  \Big/
\begin{array}{ll}
\textrm{ there are more than one different } \\ 
\textrm{minimizing geodesic connecting $x$ and $q$,}
\end{array} \right \}
\end{equation}
then, there exist a Borel set $B_q \subset \M$ with $A_q \subset B_q$ and such that for any  $q \in \M$
\begin{equation}
\label{unicidadd}
\int_{B_q} f_X(x) d \nu(x)=0.
\end{equation}
\begin{remark}
Let $C_{\M}(p)$ stand for the  cut locus of $p$ in a complete manifold $\M$, (see \cite{docarmo1992}, p. $267$).  If $P \left( X \in C_{\M}(p)\right) = 0$,  for all   $p \in \M$, then condition  (\ref{unicidadd}) is fulfilled (see \cite{pennec2006}). For instance, if the Riemannian manifold is Hadamard, then the cut locus of any point $p$ is empty and the condition is obviously fulfilled. In the unit sphere $S_{k} \defeq \{u \in \mathbb{R}^k / \Vert u \Vert= 1 \}$  the cut locus of a point $p$ is its opposite point  $-p$ and any  density defined on the sphere fulfills condition  (\ref{unicidadd}).
\end{remark}

In what follows we will assume that the Riemannian manifold $(\M,g)$ with the induced distance $d$ is connected and oriented, and that the metric space $(\M,d)$ is separable and complete. We will also assume throughout the manuscript that given two points $p,q \in \M$  there is a unique geodesic determined by $p,q$ with probability one with respect to the tensorial probability measure $d\xi(p,q) \defeq f_X(p)f_X(q)d \nu(p)d\nu(q)$.

\subsection{Geodesic balls with diameter $\overline{pq}$}
\label{bolasgeodesicas}
For any pair $p,q \in \M$ that determine a unique geodesic   $\overline{pq}$,  we define the ball of diameter  $\overline{pq}$  as the ball whose center is the middle point of the geodesic joining $p$ and $q$ with radius $d(p,q)/2$. It will be denoted by $B_{pq}$.
If $K \in \M$ is a compact set, then the family of balls 
$\mathcal{B}_{pq}$ is also a Glivenko--Cantelli class in $K$; that is, $\sup_{p,q \in K}  \vert  P \left(B_{pq} \right) - P_n\left(B_{pq} \right) \vert 
 \rightarrow 0  \quad \textrm{a.s. \ \  as  $n \rightarrow + \infty$}$, 
where  $P_n$ stands for the empirical measure of an iid sequence having distribution $P$. This last fact follows since the family of balls with center   $p$  contains the family of balls of diameter   $\overline{pq}$, which is a  Glivenko--Cantelli class, as shown in \cite{szabados1989}.
 
\section{DCOPS in Riemannian manifolds}
\label{dcops}
Let $(\M,g)$ be a Riemannian manifold and   $X: \Omega \rightarrow \M$ be a random element with density  $f_X$ with respect to $\nu$ . Given $p \in \M$, we define DCOPS at the point $p$, denoted by $BD(p)$, 
\begin{equation}
\begin{aligned}
BD(p)\defeq P\left( p \in B_{X_1 X_2} \right)  
= \int_{\M \times \M} \I_{B_{x_1 x_2}}(p)f_X(x_1) f_X(x_2) d \nu (x_1) \times d \nu(x_2),
\end{aligned}
\label{defprofundidad_man}
\end{equation}
where $X_1$ and $X_2$ are two independent random elements with the same distribution as $X$.

\begin{figure}[!ht]
\centering
\includegraphics[width=80mm]{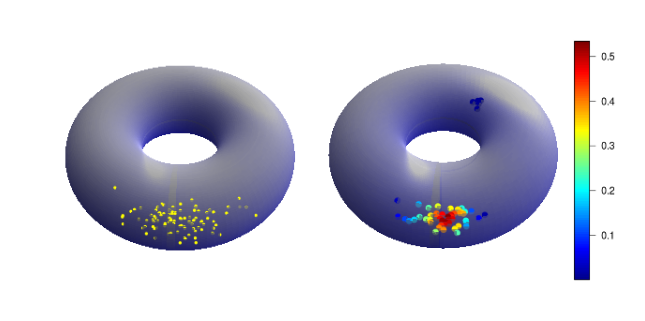}
\caption{Left--hand Panel: $100$ sample points  of the distribution  $MVM_1(\mu_1,\kappa_1, \Delta_1)$. Right--hand Panel: $\widehat{DB}$ of a sample of size $60$ from the mixture model  $0.9 MVM_1(\mu_1,\kappa_1, \Delta_1) + 0.1 MVM_2(\mu_2,\kappa_2,\Delta_2)$.}
\label{toro}  
\end{figure}

Given an iid sample $ X_1, \ldots, X_n$  having the same distribution as  $X$,  the empirical version of  $BD(p)$ is given, as before, by the order 2 associated  $U$-statistic,
\begin{equation}
\widehat{BD}_n(p)= \binom{n}{2}^{-1} \sum_{1 \leq i_1 < i_2 \leq n} \I_{B_{X_{i_1} X_{i_2}}}(p).
\label{estimadorprofundidad_man}
\end{equation}
Next, we provide two simple examples for data on the torus and on a human face. We generate data on the torus  according to the multivariate von Mises distribution (see \cite{mardia2008} and \cite{mardia2014}). The von Mises distribution, denoted by $MVM(\mu,\kappa, \Delta)$, is a distribution on the torus  $\mathbb{T}^2= S^1 \times S^1$, with density given by
\begin{equation}
f(\theta; \mu,\kappa,\Delta)= \frac{1}{Z(\kappa,\Delta)} \exp \left \{  \kappa^t c(\theta) + \frac{1}{2}s(\theta)\Delta s(\theta)  \right \}  \, (\theta \in \mathbb{T}^2).
\end{equation}
The parameters are  $\mu \in \mathbb{T}^2$ (the mean parameter), $\kappa \geq 0 \in \R^k$  (the concentration parameter) and $\Delta= (\lambda_{i,j})$ a symmetric matrix in $\R ^{k \times k}$ with  $\lambda_{i,i} = 0 \, \ \forall \ \ i=1,\ldots,k$, while $Z(\kappa,\Delta)$  is the normalization constant and $c_i(\theta)= \cos(\theta_i - \mu_i)$, $ s_i(\theta)= \sin(\theta_i - \mu_i)$, for $i=1, \ldots,k$. In the left--hand panel of the  Figure \ref{toro}, we show $100$ simulated values from a $MVM_1(\mu_1,\kappa_1, \Delta_1)$ with $\mu_1=(\pi/2,0)$, $\kappa_1=(20,20)$ and  $\Delta_1= \left( \begin{matrix}  0&1 \\ 1 &0\end{matrix}\right)$. In the right--hand panel,  we represent the empirical DCOPS values built on a sample of size $60$ from the mixture model  $0.9 MVM_1(\mu_1,\kappa_1, \Delta_1) + 0.1 MVM_2(\mu_2,\kappa_2, \Delta_2)$, with  $\mu_2=(7/4 \pi,0)$, $\kappa_2=(100,100)$ and  $\Delta_2= \Delta_1$. We can see how the depth decreases from the center $\mu_1$ and the small depth values for the sample of the second component of the mixture (outliers).

\begin{figure}[!ht]
\centering
\subfloat{\includegraphics[width=54mm]{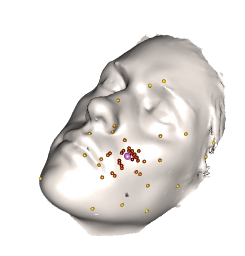}}
\subfloat{\includegraphics[width=50mm]{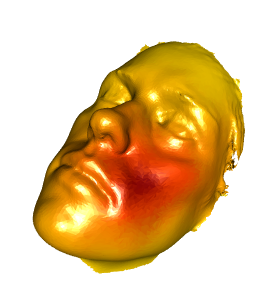}}
\caption{Left--hand panel: $\widehat{DB}$ at the 50 sample points. In violet the deepest point. Right--hand Panel: Image of  the intensity depth to the deepest point $p$ according to $\widehat{DB}$. A larger intensity of red indicates a larger depth.} \label{cara1}  
\end{figure}

We can also consider a manifold to be a surface, such as a human face. Given a fixed point $p$ (represented in violet in the left--hand panel of Figure \ref{cara1}) on the human face $\mathcal{M}$ and $d$ the geodesic distance, we define a density $f$ on $\mathcal{M} $ by
\begin{equation}
\label{densidad_cara}
f(x) \defeq \frac{\psi(x)}{ \int_{\mathcal{M}}  \psi(t) d \nu(t)},
\end{equation}
with $\psi(x)= \frac{1}{1+d(p,x)}$. We chose a sample of $50$ points at random with respect to a discretization of (\ref{densidad_cara}) in a grid of points on the face. We calculate the depth of these points (see Figure \ref{cara1}, left--hand panel) as well as the depth of all of the points in the face with respect to the sample depth (see Figure \ref{cara1}, right--hand panel).

\subsection{Some properties of DCOPS on Riemannian manifolds}
\label{pdcops}
\begin{figure}[!ht]
\centering
\includegraphics[width=100mm]{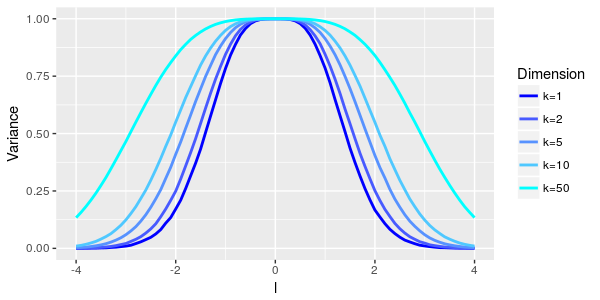}
\caption{Graphs of  marginal variances $\sigma_k^2(lu)=4  P_2(D_{lu}) \left( 1- P_2(D_{lu})\right)$ functions with respect to $l\in [-4,4]$ for  $k= 1,2,5,10 \,\, \textrm{and} \, \, 50. \, (u \in \mathbb{R}^k,\, \Vert u \Vert=1)$ and $D_{lu} \defeq \{ (x_1,x_2) \in \mathbb{R}^k \times \mathbb{R}^k: B_{x_1 x_2} \ni lu \}$}
\label{varianzas}  
\end{figure}

\begin{enumerate}
\item (\textit{Orthogonal invariance.})
DCOPS is invariant with respect to orthogonal transformations in the space  $\R^k$ where it is embedded.
\item (\textit{Vanishing at infinity.})
If we fix a point $p$ in the manifold $\mathcal M$ and the distance of another point $q$ diverges to infinity, then $BD(q)$ tends to zero.

This last result is shown in the following theorem.
\begin{theorem}
Let  $q$ be a fixed point in $\M$. Then, we have that 
\label{desvanecimientoinfinitomanifold}
\begin{enumerate}
\item[(a)] $\sup_{d(p,q)>M} BD(p) \rightarrow 0 $   when  $M \rightarrow +\infty$  
\item[(b)] $ \lim_{M\to \infty} \sup_{d(p,q)>M} \widehat{BD}_n(p)=0$ a.s.
\end{enumerate}
\end{theorem}
The proof is given in the Appendix.

\item (\textit{Continuity})
The next theorem shows that if two points are close with respect to the geodesic distance, then their respective depths are close. 

Given two points $p, q \in \mathcal M$ we denote by $S_{pq}$  the geodesic sphere of diameter $\overline{pq}$. That is,
 $S_{pq} \defeq \{x \in \M \ \ \mbox{such that} \ \  d(\frac{pq}{2},x)= d(p,q) /2\}$.

We will assume the following condition.

\bf Condition VC. \rm Let $\mathcal M$ be a compact manifold. $\mathcal M$ satisfies the  \textbf{VC} condition if the family of sets $\mathcal{A}=\{ A_p: p \in \M \}$, with 
 $$A_p= \left \{(x,y) \in \M \times \M \ \ \mbox{such that} \ \  \frac{xy}{2} \in \bar{B}\left(p, \frac{d(x,y)}{2}\right) \right \},$$ has a finite  Vapnik--Chervonenkis dimension ($VC$--dimension),  see \cite{steele1975} and  \cite{dudley1978}.
\begin{theorem}
\label{teoremacontinuidadmanifold}
Given $p,q \in \M$ we have that,
\begin{enumerate}
\item [(a)] $\vert BD(p) -  BD(q) \vert \rightarrow 0  \quad \textrm{as} \quad d(p,q) \rightarrow 0$. 
\item [(b)] Let $\M$ be a compact manifold fulfilling condition \textbf{VC}. Assume that $X$ has a bounded density $f_X$. Then, we have that for $\epsilon>0$
\begin{displaymath}
 \sup_{d(p,q)<\epsilon} \vert \widehat{BD}_n(p) - \widehat{BD}_n(q)  \vert  < \gamma(\epsilon)  + R_n,
\end{displaymath}
 where $\gamma(\epsilon)$ is a deterministic function that converges to $0$ when $\epsilon \rightarrow 0$ and $R_n$ converges a.s. to $0$ when  $n \rightarrow \infty$.
\end{enumerate}
\end{theorem}
\begin{remark}
\label{vc2}
 In the Euclidean case,  $\mathbb{R}^k$ the condition {\bf {VC} } holds and compactness is not necessary. Indeed, the $VC$--dimension $\nu_{\mathcal A}$ is bounded by $3k+1$. As a matter of fact
, since
\begin{align*}
A_p= \{(x,y) \in \mathbb{R}^k \times \mathbb{R}^k /  \langle p,x+y \rangle - 2 \langle x,y\rangle+ 2\Vert p \Vert^2 \geq 0  \} \\
&\hspace{-8 cm}= \{(x,y) \in \mathbb{R}^k \times \mathbb{R}^k / g_p(x,y)  \geq 0  \}
\end{align*}
with $g_p(x,y)=\langle p,x+y \rangle - 2 \langle x,y \rangle+ 2\Vert p \Vert^2$. But,  $\mathcal{G}=\{g_p\}_{p \in \mathbb{R}^k}$ is the linear span generated by the set of functions, $g_0(x,y)= 1$;   $g_i(x,y)= x_i  \quad  i=1, \ldots,k$ ;  $g_j(x,y)= y_j  \quad  j=1, \ldots,k$;  $g_t(x,y)= x_t y_t  \quad  t=1, \ldots,k$.
and the results in \cite{steele1975} and  \cite{dudley1978} imply that  $\mathcal{V}_\mathcal{A}  \leq 3k+1$. 
\end{remark}
The case of FDA is analyzed in subsection \ref{fda} below.

\begin{remark}
 Research in the area where the \textbf{VC} condition holds for data on manifolds is in its early stages, and a first approach can be found in \cite{ferri2008}. For some simple manifolds like the sphere, the torus and the cylinder it can be proved directly as follows.  Let $S^2$ the sphere embedded on $\mathbb{R}^3$ with center $(0,0,0)$ and radius $1$. Let $x=(x_1,x_2,x_3), y=(y_1,y_2,y_3)$, then $(x,y) \in A_p$ if it is  satisfied 
$$g_p(x,y) \defeq  \Vert \frac{x-y}{2}\Vert^2 - \sum_{i=1}^3 \left( p_i -\frac{x_i+y_i}{2} \right)^2 >0,$$ with $p=(p_1,p_2,p_3)$. Therefore $\mathcal{G}$ is generated by  $g_0(x,y)= 1; g_1(x,y)= \Vert \frac{x-y}{2}\Vert^2; \ g^{(k)}_i(x,y)= \left(\frac{x_i+y_i}{2} \right)^k, i=1,2,3; k=1,2.$
For the torus case,  let $x=(x^{(1)},x^{(2)}), y=(y^{(1)},y^{(2)}) \in  S^1 \times S^1$,   $(x,y) \in A_p$ if it satisfied 
$$g_p(x,y) \defeq  d^2(x,y)/4 - d^2_{S^1} \left(p^{(1)}, \left(\frac{xy}{2}\right)^{(1)} \right) -  d^2_{S^1} \left(p^{(2)}, \left(\frac{xy}{2} \right)^{(2)} \right) >0.$$
Recall that the relationship between the chord and geodesic distance is
$$d_{S^1}(x^{(k)},y^{(k)})= K d_{\mathbb{R}^2}(x^{(k)},y^{(k)})$$
 where $K$ a known constant. Hence, a similar argument as the one used for  the sphere, shows that  $\mathcal{G}$ is of finite $VC$--dimension on the torus. The case of the cilinder is completely analogous.
\end{remark}
\item   (\textit{Uniform convergence.}) We start by proving the uniform convergence of the empirical version of the depth to the population version. For a fixed
 $p \in \M$, the estimator $\widehat{BD}_n(p)$  proposed in (\ref{estimadorprofundidad_man}) is unbiased and the strong pointwise consistency follows from the fact that the kernel   $h(x,y)=\I_{B(x,y)}(p)$ of the  $U$-statistic is bounded (see \cite{serfling1980}, Theorem A, [5.6]). We will show  the uniform convergence in the next theorem. The proof is given in the Appendix. 
\begin{theorem}
\label{consistenciamanifold}
Let $\M$ be a compact manifold fulfilling condition \textbf{VC}. Then, we have
$\sup_{p \in \M}  \vert \widehat{BD}_n(p) - BD(p) \vert \rightarrow 0, \quad \textrm{a.s.} \quad \textrm{for} \quad n \rightarrow + \infty$
\end{theorem}
\end{enumerate}
The asymptotic distribution of the depth function $\widehat{BD_n}$ can  be derived using the results  in \cite{arcones1993} for $U$-statistics. 
The set $\mathcal{D}$ is defined in (\ref{D}). Let $\ell^\infty (\mathcal{D})$ be the set of all bounded functions  $f: \mathcal{D} \to \mathbb{R}$. 
We consider the empirical depth $U$-process indexed by $\mathcal{D}$, $\widehat{BD}_n \defeq \{U_2^n\left(D_x \right): D_x \in \mathcal{D} \}$ , that is, for each $x$, $U_2^n\left(D_x \right)$ is a $U$-statistics.
\begin{theorem} [Asymptotic distribution] 
\label{normalidad}
Let $\mathcal M$ be a compact manifold fulfilling condition \textbf{VC}. Let  $ \{ X_n : n \geq 1 \}$ be an iid sequence distributed as  $X$. Let $P$ be the distribution of $X$. Then, the stochastic process $n^{1/2} \left(\widehat{BD}_n -BD \right)$ converges in law  to  $2G_P \, \textrm{in} \ \  \ell^\infty(\mathcal{D})$, where $G_P$ is the  Brownian bridge associated with  $P$ and indexed by $\mathcal{D}$. That is, $G_P$ is a centered Gaussian process  with covariance $E \left[ G_P(D_x) G_P(D_y)\right]= P_2(D_x \cap D_y)- P_2(D_x)P_2(D_y), \, (D_x,D_y \in \mathcal{D})$,
where $$P_2(D) \defeq \int_{\mathbb{R}^k \times \mathbb{R}^k} \I_{D}(x_1,x_2) \, dP(x_1) \times dP(x_2), \,\, (D \in \mathcal{D}).$$
\end{theorem}
Notice that, the convergence in law involved in  Theorem (\ref{normalidad}) is in the sense of  \cite{hoffmann1991}.
\\
In particular, for a fixed $x$ different from the deepest point, the marginal distribution of the limit process $2G_P$ is given by $2G_p(D_x) \sim N \left[ 0, \sigma_k^2(x) \right]$ with $\sigma_k^2(x)=4  P_2(D_x) \left( 1- P_2(D_x)\right)$. 
\\
In the  Euclidean space $\mathbb{R}^k$,  if $X \sim N(0, I_k)$ and we denote by $S(x,y)= \{ t \in \mathbb{R}^k : \langle y-x,x-t \rangle \geq 0\}$, 
$$ P_2(D_x)= \frac{1}{(2\pi)^{k/2}} \int_{\mathbb{R}^k} g_x(y) e^{- \Vert y\Vert^2 /2}dy, $$
where  $$g_x(y)=P \left( S(x,y) \right)= \frac{1}{\sqrt{2\pi}} \int^{+ \infty}_{ \frac{\langle x-y,x \rangle} {\Vert x-y \Vert}} e^{-t^2/2}dt.$$
We display in the Figure \ref{varianzas} the marginal variance  $\sigma_k^2(lu)$ as a function of $l \in \mathbb{R}$. Here, $u$  is any unitary vector of $\mathbb{R}^k$. Notice that the variance increases with the dimension of the space.
If the distribution is spherical symmetric around  $c$ (the deepest point) the marginal distribution at $c$ corresponds to the degenerate case of a $U$-statistic and $n\left(\widehat{BD}_n(c) -BD(c) \right)$ converge in distribution to $\sum_{j=1}^{+\infty} \lambda_j (\chi^2_{j} -1)$ where $\chi^2_{j}, (j \geq 1)$ are independent chi--square random variables with one degree of freedom (see \cite{serfling2000}).
  
\subsection{\textbf{A consistency result for FDA}} \label{fda}
In a separable Hilbert space $\mathcal{H}$, as we will show in  Theorem  \ref{hilb} below,  the \textbf{VC} condition is not necessary to prove the uniform consistency. Instead, we will need the following weaker assumption:

\textbf{HB) B-continuity} A probability measure $P$ defined on a separable Hilbert space $\mathcal{H}$ fulfils \bf HB \rm if $P(\partial B(y,r))=0$ (Here $\partial A$ denotes the boundary of the set $A$) for all $r>0$ and $y\in \mathcal{H}$.


\begin{theorem}
\label{hilb}
 Let $\{X_n\}_n$ be a sequence of random elements with common distribution $P_n$, valued in a separable Hilbert space $(\mathcal{H},\|\cdot\|)$.  Let $P$ be a probability distribution fulfilling \textbf{HB}. Then, if   $P_n\rightarrow P$ weakly,
$$\sup_{p \in \mathcal{H}}  \vert \widehat{BD}_n(p) - BD(p) \vert \rightarrow 0, \quad \textrm{a.s.} \quad \textrm{for} \quad n \rightarrow + \infty$$
\end{theorem} 
In order to prove Theorem  \ref{hilb} we will use the following  result proved in \cite{billingsley1967} (it still holds when $\mathcal{H}$ is a separable Banach space), see Theorem 1 and Example 3  in \cite{billingsley1967}. \\

\begin{theorem}
\textbf{(Billingsley and Tops{\o}e).}  \textit{Let $S$ be a separable metric space and $\mathcal{B}(S)$ the class of all bounded, real,  measurable functions defined on $S$. Let $\mathcal{F} \subset \mathcal{B}(S)$ a subclass of functions, then 
$$
 \sup_{f \in \mathcal{F} } \left \vert \int f dP_n  - \int f dP \right \vert  \rightarrow 0,
$$ for every sequence $P_n$ that converge weakly to $P$, if and only if
\begin{equation}\label{BT1} 
\sup \left \{ \vert f(z)-f(t) \vert : f \in \mathcal{F},z,t \in S \right \} < \infty,
\end{equation} and for all $\epsilon >0$, 
\begin{equation} \label{BT2}\lim_{\delta \rightarrow 0} \sup_{f \in \mathcal{F}} P\left( \{ x: \omega_f \left( B(x,\delta) \right) > \epsilon \} \right)=0
\end{equation} where $\omega_f(A)=\sup \left \{ \vert f(x) -f(y )\vert: x,y \in A \subset S \right \}$ and $B(x, \delta)$ is the open ball of radius $\delta>0$.}
\end{theorem}

The following two corollaries follows inmediately.
\begin{corollary} \label{cor00} Consistency of the deepest point. Under the hypotheses of Theorem \ref{hilb}, if $BD(p)$ has a unique minimum under $P$,  then 
$$ \hat\theta_n:= argmax_{p} \widehat BD_n(p)  \,\, \textrm{ converge to } \,\, \theta:=argmax_{p}  BD(p) \ \mbox {a.s.}.$$  
\end{corollary}

 \begin{corollary}  Robustness of the deepest point. Under the hypotheses of Corollary \ref{cor00}, $\hat{\theta_n}$ is qualitatively robust in the sense of \cite{hampel1971}.
\end{corollary}

\section{Some simulated examples}
\label{simulated}
We provide a small simulation study for data in, the sphere, and data on the cone of positive definite matrices (manifold). We will end this section by considering some examples in the Euclidean space $\mathbb R^k$ with $k=5$ and $k=20$. These last simulations will illustrate that the DCOPS approach of depth is versatile and can be performed even for large dimension.
\subsection{Simulated examples for data on the sphere}
\begin{figure}[!ht] 
\centering
\subfloat{\includegraphics[width=56mm]{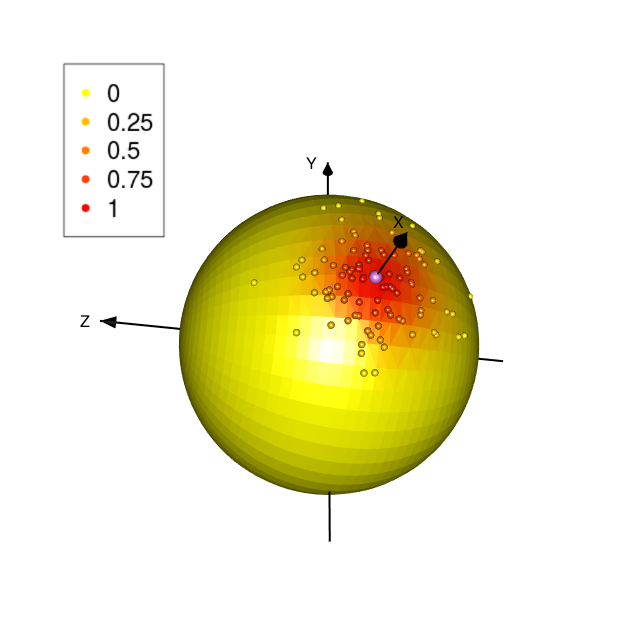}}
\subfloat{\includegraphics[width=56mm]{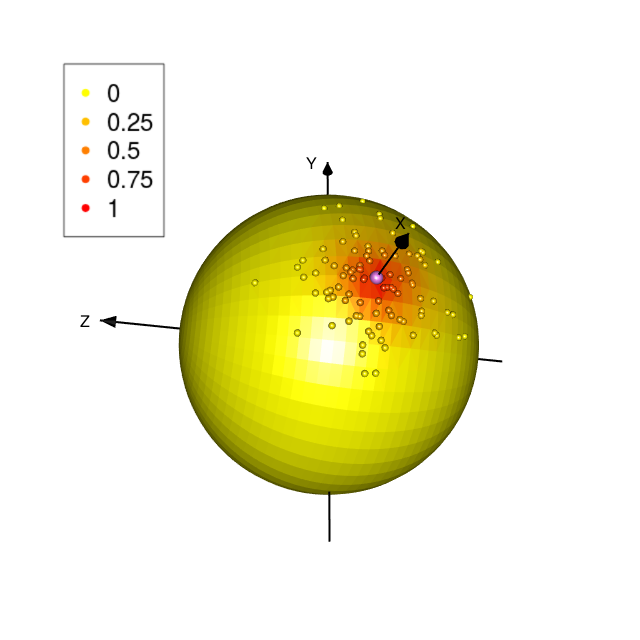}}
\caption{ Empirical depths for a sample of size $100$ for data with  Von Mises--Fisher distribution with $\kappa =15$  and $\mu= (1,0,0)$. Higher intensities of red indicate larger depth values.  Left--hand panel:  Normalized DCOPS. Right--hand panel:  Normalized ATD.} \label{prof_esfera}  
\end{figure}

\begin{figure}[!ht]
\centering
\subfloat{\includegraphics[width=58mm]{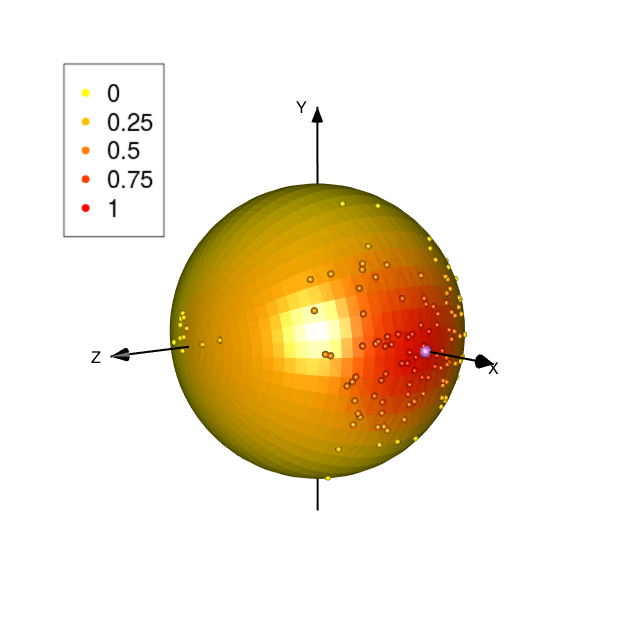}}
\subfloat{\includegraphics[width=58mm]{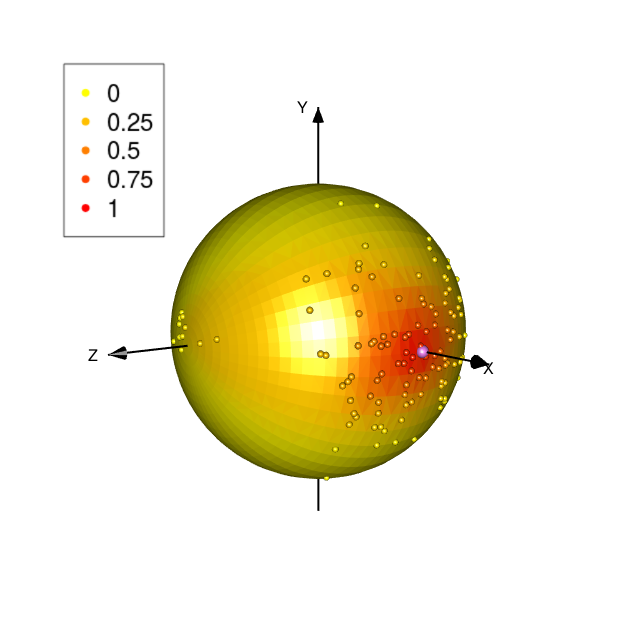}}
\caption{ Empirical depths for the mixture distribution with sample size $120$. Left--hand panel: DCOPS. Right-hand panel:  Normalized ATD. Higher intensities of red indicate larger depth values.}  \label{prof_esfera_out}  
\end{figure}

We consider a sample of iid data of size $100$ on the sphere $S^2=\{ x \in \mathbb{R}^3 /  \Vert x \Vert= 1\}$ with Von Mises--Fisher distribution given by
$f(x,\mu,\kappa)= C(x) e^{\kappa \mu^t x} \I_{ S^2}(x)$,
where $\kappa \geq 0$  and $\mu \in S^2$ are the concentration and the directional media parameters, respectively. $C(x)$ is the normalizing constant, see \cite{mardia2000}.
In Figure \ref{prof_esfera} we plot the depth at  $100$ values generated with $\kappa =15$  and $\mu= (1,0,0)$.  We observe larger values of the depth (red color) at a neighborhood of the directional mean, and it decreases (yellow color) at the antipodal point for both DCOPS and the extension of Tukey depth for data on the sphere ATD (angular Tukey's depth) introduced in \cite{liu1992b}. Figure \ref{prof_esfera} shows the similarity of both depth measurements (normalized between 0 and 1) at the sample points. However, the computational time to compute  DCOPS is considerably smaller.
Next, we add outliers considering a mixture of von Mises--Fisher distributions given by  
$f(x)= 0.9 f(x,(1,0,0),10) + 0.1 f(x,(0,0,1),50)$.
In Figure \ref{prof_esfera_out}, we perform the same plot for a sample of 120 points and we can see that both depths give small values to the outliers, which are on the left--hand hemisphere. The behaviour of both depth functions is quite similar in this case.

\subsection{Simulated examples for positive definite  random matrices}
Several statistical problems  involve the study of data that are a collection of positive definite matrices. For example, the robust method estimation  for covariance matrix or principal components analysis (see \cite{chen2017} and \cite{ggordaliza2017}), learning and shape analysis, among others.  The space of such matrices is a Riemannian manifold that is denoted by  $(\mathbb{P}_k,g)$. Given two matrices   $A,B \in \mathbb{P}_k$, there exist a unique geodesic determined by $A$ and $B$ (see \cite{moakher2005}), which is given by,
\begin{equation}
\gamma(s) \defeq  A ^{1/2} \left(A ^{-1/2} B A ^{-1/2}  \right) ^ s A ^{1/2}, 
\end{equation}  
and we can find the middle point between $A$ and $B$ that will be denoted by
 $A\#B$, as well as the geodesic distance among them,
\begin{align*}
A\#B = A ^{1/2} \left(A ^{-1/2} B A ^{-1/2}  \right) ^{1/2} A ^{1/2}, \\
&\hspace{-6.2 cm}d(A,B)= \Vert \log  \left(A ^{-1/2} B A ^{-1/2}  \right) \Vert,
\end{align*}
where  $\Vert \cdot \Vert$ is the Hilbert--Schmidt norm. We observe that the middle point $A\#B$ is nothing but the geometric mean between the matrices. 

We start by considering the data generated by a Wishart distribution $W_3(\Sigma, m)$ in  $\mathbb{P}_3$, with parameters $m=20$ and $\Sigma=I_3$. As is well known, a random matrix with Wishart distribution $S$ with parameters $\Sigma$ and $m$ can be generated from iid multivariate normal vectors $N(0,\Sigma)$ by $S \defeq \sum_{i=1}^m X_i X_i^t$. Its expectation is given by $\mathbb E(S)=m \times \Sigma$. 

\begin{figure}[!ht] 
\centering
\subfloat{\includegraphics[width=45 mm]{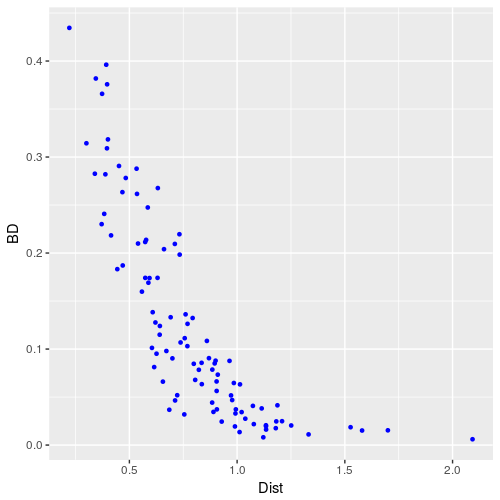}}
\subfloat{\includegraphics[width=68 mm]{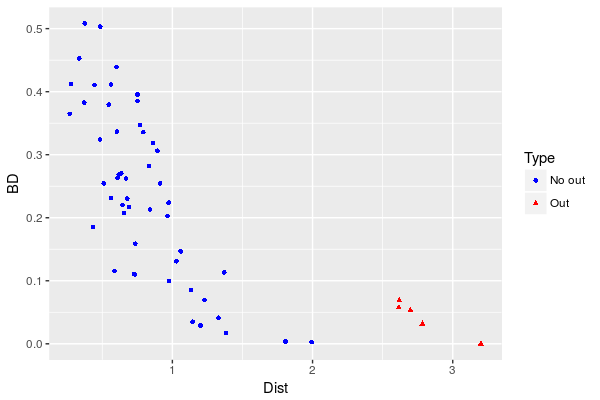}}
\caption{  Left--hand panel: Plot of $\widehat{BD}_n$ as a function of the geodesic distance to the expected value for a sample of size 100 from a $W_3(I_3,20)$ distribution. Right--hand panel: Plot of $\widehat{BD}_n$ as a function of the geodesic distance to the expected value of the non--contaminated data. The outliers are marked in red.} \label{prof_matrices1} 
\end{figure}

We generate a sample of size $100$. For each matrix in the sample, we calculate the geodesic distance to its expectation and its DCOPS depth $\widehat{BD}_n$.  The results are given in Figure \ref{prof_matrices1}. It can be observed that the depth decreases as we go further from the expected value.  Next, we consider a mixture of two Wishart distributions 
$ 0.9  W_3(I_3,20) + 0.1  W_3(I_3/10,50)$, where the outliers are generated from the distribution $  W_3(I_3/10,50)$ and a sample size of $120$.
Figure \ref{prof_matrices1} shows that $\widehat{BD}_n$ gives small values to the outliers. 
\newpage
\subsection{Simulated examples in $\mathbb{R}^k$ with $k=5$ and $k=20$}
\begin{figure}[!ht]
\centering
\subfloat{\includegraphics[width=68mm]{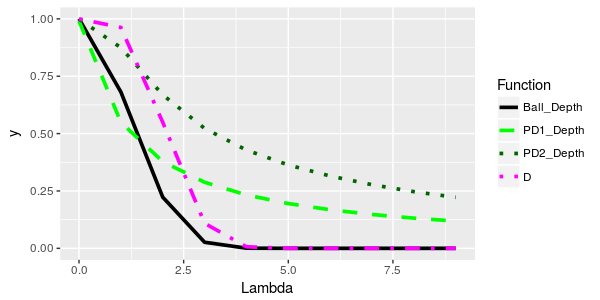}}
\subfloat{\includegraphics[width=68mm]{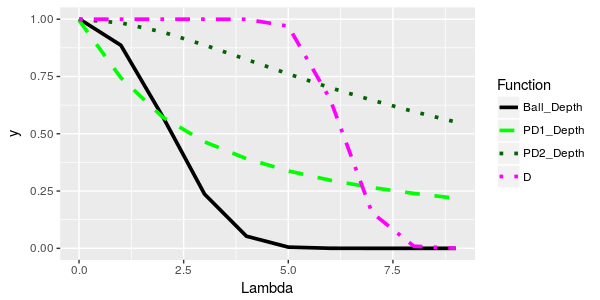}}
\caption{Plot of normalized depth $PD_{1,n},4PD_{2,n}, 2BD_n$ 
and  the probability $D$ for $N(0,I_k)$ 
data restricted to the half--line   $y=\lambda (1,0,0, \ldots,0) \in \mathbb{R}^k$, $\lambda>0$ as a function of the distance from the origin. Left--hand Panel : $k=5$ and $n=1000$.  Right--hand Panel : $k=20$ and $n=1000$ }\label{prof_alta}  
\end{figure}

As mentioned before for large $k$, empirical simplicial and half--space depths are computationally hard to compute. Hence, we compare DCOPS with  two other notions of depth that has been introduced in \cite{liu1992} and \cite{fraiman2009}, respectively. These notions are based on one--dimensional projections.
\\
Let $\X_1, \ldots, \X_n \in \mathbb R^k$ be iid random vectors with the same distribution as $\X$.   The first empirical depth notion is defined as follows. 
For each $\textbf{u} \in \mathbb R^k$ ($\Vert \textbf{u} \Vert =1$), we project the sample onto the direction $\textbf{u}$ and obtain a one dimensional sample
$\aleph_{n,u} = \{X_{i,u}= \langle \X_i, \textbf{u}\rangle, \ \   i=1,\ldots,n\}$.
A measure of ``outlieness" is then defined by 
\begin{equation}{\label{ou}} OU_n(\mathbf x)  \defeq \sup_{\substack{\Vert \textbf{u} \Vert=1}}   \frac{\vert \langle \mathbf x,\textbf{u} \rangle  - \mu_{\aleph_{n,u}}  \vert}{\tau_{\aleph_{n,u}}}, \ \ \mathbf x \in \mathbb R^k,
\end{equation}
where $\mu_{\aleph_{n,u}}$ and $\tau_{\aleph_{n,u}}$ are the median  and  the median absolute deviation  of the projected sample $\aleph_{n,u}$ respectively (see \cite{maronna2006}, page $5$). The corresponding population version is obtained by replacing  $\mu_{\aleph_{n,u}}$ and $\tau_{\aleph_{n,u}}$ by $\mu_u$ and $\tau_u$ the median and the absolute deviation of the random variable $\langle X, \mathbf u\rangle$ in (\ref{ou}). 
In \cite{liu1992} the author suggests to use as  empirical depth measure 
$PD_{1,n} (x) \defeq  \left( 1+OU_n(x) \right)^{-1}$.
In practice, $PD_{1,n} (x)$ can be approximated by taking the maximum over a large sample of random directions instead of the supremum over the whole unit sphere.  
\\
The second depth notion introduced for the functional case in \cite{fraiman2009} can also be used for high dimensional data. 
\\
If we denote by  $g_{\textbf{u}}(x)=\langle x, \textbf{u} \rangle$, then the depth is defined by 
 \begin{equation}
PD_{2}(x) \defeq   \int_{S_k} F_{g_{\textbf{u}}(X)} \left(g_{\textbf{u}}(x)\right)  \left(   1-  F_{g_{\textbf{u}}(X)} \left(g_{\textbf{u}}(x)\right)  \right) d\textbf{u},
 \end{equation}
where $F_{g_{\textbf{u}}(X)} (\cdot)$ stands for the cumulative distribution function of the random variable $g_{\textbf{u}}(X)$, which can be  estimated  as the average for $K$ random directions (uniform on the unit sphere); that is, 
$$PD_{2,n}(x) \defeq   \frac{1}{K} \sum_{i=1}^K \hat{F}_{g_{\textbf{u}_i}} (g_{\textbf{u}_i}(x)  \left(1-  \hat{F}_{g_{\textbf{u}_i}} (g_{\textbf{u}_i}(x))\right),$$
where $\hat{F}_{g_{\textbf{u}_i}}(\cdot)$ is the empirical distribution function of the random sample  $$\{g_{\textbf{u}_i}(X_1), \ldots, g_{\textbf{u}_i}(X_n) \}.$$
\\
First, we consider the sample sizes  $n=1000$ in dimensions $k=5,20$, for $N(0, I_k)$ data. 
Figure  \ref{prof_alta} depicts the depth of the points lying on the half line  $y=\lambda (1,0,0, \ldots,0) \in \mathbb{R}^k$ with $\lambda >0$ for the normalized depth $PD_{1,n}$, $4PD_{2,n}$ and $2BD_n$  with respect to the Euclidean distance to the origin. The same plot also depicts the function  $D(\lambda)= P(\Vert X \Vert >\lambda), \, \left( X \sim N(0,I_k) \right)$. In both cases a good behavior of DCOPS is observed, the depth decreases when the point moves away from the origin and quickly converges to zero.

\subsubsection{Contaminated data}
\begin{figure}[!ht]
\centering
\subfloat{\includegraphics[width=100mm]{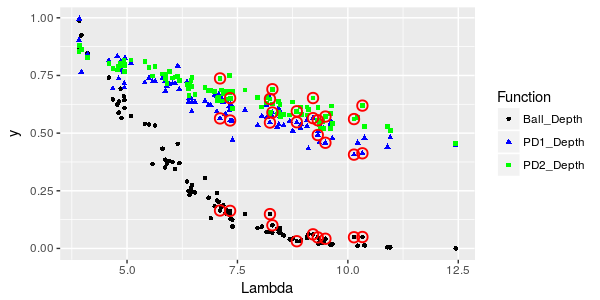}}
\caption{ Plot of the empirical normalized depths $PD_{1,n}$, $4PD_{2,n}$, $2BD_n$  at $100$ points generated at random as a function of the distance to the origin. The outliers are marked in red.} \label{prof_alta_out}  
\end{figure}

We now consider a contaminated sample given by 
\begin{equation}
 (1-\alpha) N_2\left(\mu_1, \Sigma_1 \right) +  \alpha N_2\left(\mu_2, \Sigma_2\right),
\label{contaminada}
\end{equation} 
A sample of size $n=5000$ was generated with a mixture of multivariate normals given by equation (\ref{contaminada}) with
 $\mu_1=(0,0, \ldots,0), \mu_2= 2 (1,1, \ldots,1)  \in \mathbb{R}^{10}$,  $\Sigma_1= \Sigma_2= I_d$ and $\alpha=0.1$. In Figure  \ref{prof_alta_out}, we plot the empirical normalized depths  $PD_{1,n}$, $4PD_{2,n}$ and $2BD_n$ having at hand a set of $100$ points generated at random with the same contaminated  distribution 
as a function of the distance $\lambda$ to the origin.   The outliers are marked in red. It can be seen that all three depths give moderate or small values to the outliers, although  $BD_n$ seems to outperforms the competitors.
\section{Concluding remarks}
We propose a notion of depth for data in Riemannian manifolds such that:
\begin{itemize}
\item It shares the classical desirable properties generally assumed in finite dimensional spaces.
\item Its definition as well as its  main properties hold in the frame of a separable Hilbert space. Strong uniform convergence is shown also for FDA.
\item We have shown the strong consistency of the empirical depth. We also  provide the asymptotic distribution for the empirical estimator of this  depth that covers the important case where the data live in a Riemannian manifold.
\item We provide several simulated examples in Riemannian manifolds frames (such as the sphere, the torus, and the space of positive definite matrices). 
\item We show here, through simulations,  that it can be calculated easily  for high dimensional spaces (indeed its computational complexity is of order $k \times n^2$, $n$ being the sample size), maintaining a good performance. In general, the computational times grow in a Riemannian manifold  due to the geodesic calculation. However, we consider that it is the methodology appropriated in this paradigm.

\item It is important to note that other very  interesting  depth functions have been proposed in more general spaces (separable metric spaces), see for example  \cite{liu2011} and \cite{carrizosa1996}. However, after defining the concept of statistical depth in the general space, the desirable properties are  often not studied. We postpone to a future work the study of the consistency and the CLT of these estimators in metric spaces. Further, we will also compare through simulations the performance of these depths for Riemannian manifolds with respect to DCOPS. Some recent comparative simulation study in $\mathbb{R}^k$ might be found in  \cite{shahsavarifar2018}.

\end{itemize}
\section*{Appendix}
\begin{enumerate}[label=(\Alph*)]

\item \textbf{Proof of Theorem \ref{desvanecimientoinfinitomanifold}.}
\begin{enumerate}[label=(\alph*)]
\item We start by showing that if  $d(p,q)>M$, 
\begin{equation}
\label{des1}
BD(p)=P\left( p \in B_{X_1 X_2} \right) \leq 2 P\big( d(q,X)>M/ 4 \big),
\end{equation}
which will be shown by contradiction. Given $x_1, x_2 \in  B(q, M/4)$, let us denote by $\frac{x_1x_2}{2}$  the middle point of the geodesic determined by  $x_1$ and $x_2$. Then,  
\begin{displaymath}
d(q,\frac{x_1x_2}{2}) \leq d(q,x_1) + d(x_1,\frac{x_1x_2}{2}) = d(q,x_1) + \frac{d(x_1,x_2)}{2} < \frac{M}{2}.
\end{displaymath}
Note that then $p \notin B_{x_1 x_2}$, since if  $d(p,\frac{x_1x_2}{2})<\frac{d(x_1,x_2)}{2}$ then,
\begin{displaymath}
d(q,p) \leq  d(q,\frac{x_1x_2}{2})+d(\frac{x_1x_2}{2},p)< \frac{M}{2} + \frac{M}{4}< M,
\end{displaymath}
which is a contradiction. Since the space is separable, $d(X,q)$ is a random variable and, therefore,

\begin{align*}
BD(x)\leq P \left( \left \{ X_1 \notin B(q,M/4) \right \}\cup \left \{ X_2 \notin B(q,M/4) \right \} \right) \\
& \hspace{-7.2 cm}\leq  2 P\left( d(X,q) > M/ 4 \right).
\end{align*}

Finally, using the dominated convergence Theorem, we get that 
$$\lim_{M \to \infty} P\left( d(X,q) > M/ 4 \right) =0. \ \ \QEDB $$ 
\item Now suppose that for any $M>0$,there exists $\Omega_1 \subset \Omega$ of positive probability $\delta=P( \Omega_1)>0$, for which there exist $ \gamma>0$  and $ x_0 \in \M$ with $d(q,x_0)>M$ such that $\widehat{BD}_n(x_0,w) > \gamma>0$ for all $w \in \Omega_1$. 
Choose $M_0$ such that $P\left(d(q,X)> M_0/4 \right)< \frac{\delta}{2n} $.
 Since $\Omega_1 \subset \bigcup_{i=1}^n \left\{ d(q,X_i) > M_0/4 \right\}$ we get $\delta \leq  n P \left\{ d(q,X) > M_0/4 \right\} <\delta/2$, which is a contradiction. 
\QEDB
\end{enumerate}
\item \textbf{Proof of Theorem \ref{teoremacontinuidadmanifold}.}
\begin{enumerate}[label=(\alph*)]
\item Under our assumptions, we have that $P(p \in S_{X_1X_2})=0$
 for all $p \in \M$. Now consider the difference,
 \begin{align*}
 \vert BD(p) -  BD(q) \vert = \left \vert P \left( p \in B_{X_1X_2}  \right) - P \left( q \in B_{X_1X_2}  \right)\right \vert \\
 & \hspace{- 9.5 cm}= \left \vert P \left( p \in B_{X_1X_2}, q \notin B_{X_1X_2}  \right) - P \left( q \in B_{X_1X_2}, p \notin B_{X_1X_2}  \right) \right \vert \\
 & \hspace{- 9.5 cm} \leq \max \left \{ P \left( p \in B_{X_1X_2}, q \notin B_{X_1X_2}  \right) ,  P \left( q \in B_{X_1X_2}, p \notin B_{X_1X_2}  \right) \right\}.
 \end{align*}
The event that the ball $B_{X_1X_2}$ contains only one of the points $p, q$ is contained in the event that the sphere $S_{X_1 X_2}$  intersects the geodesic  $\overline{pq}$ determined by $p$ and $q$; hence, from the dominated convergence theorem we get 
$$
\vert BD(p) -  BD(q) \vert \leq P \left(S_{X_1X_2} \cap \overline{pq} \neq \emptyset \right) \rightarrow 0 \ \mbox{as} \ d(p,q) \rightarrow 0.  \ \ \QEDB
$$ 
\item  Let $R^{*}= diam(\M)$ and denote by $A \triangle B$ the symmetric difference between the sets. Then, 
\begin{align*}
\vert \widehat{BD}_n(p) - \widehat{BD}_n(q)  \vert \leq  \binom{n}{2}^{-1} \sum_{1 \leq i_1 < i_2 \leq n} \left \vert \I_{B_{X_{i_1} X_{i_2}}}(p)-\I_{B_{X_{i_1} X_{i_2}}}(q) \right \vert \\
&\hspace{-12cm} \leq  \binom{n}{2}^{-1} \sum_{1 \leq i_1 < i_2 \leq n} \left \vert \I_{B\left( p,d(X_{i_1},X_{i_2})/2 \right) \triangle B\left( q,d(X_{i_1},X_{i_2})/2 \right)  } \left( \frac{X_{i_1}X_{i_2}}{2} \right) \right \vert.
\end{align*}

Now observe that  $$h_{p,q}(x,y)=  \I_{B\left(p,d(x,y)/2\right) \triangle B\left( q,d(x,y)/2 \right)  } \left( \frac{xy}{2} \right)$$ is a kernel of a $U$-statistic of order $2$.  We can bound  $\vert \widehat{BD}_n(p) - \widehat{BD}_n(q)  \vert$ by 
\begin{small}
\begin{align*}
\binom{n}{2}^{-1} \sum_{1 \leq i_1 < i_2 \leq n} \left \vert h_{p,q}\left(X_{i_1},X_{i_2}\right) - E \left[ h_{p,q}\left(X_{i_1},X_{i_2}\right) \right] \right\vert + E \left[ h_{p,q}\left(X_{1},X_{2}\right) \right].
\end{align*}
\end{small}
On the other hand, since the family of sets $$\mathcal{A}_2=\{ A_{pq}=A_p \triangle A_q \}_{(p,q) \in E\times E}$$  has a finite $VC$--dimension, Corollary 3.3 in   \cite{arcones1993} entails that 
\begin{align*}
R_n \defeq \\
& \hspace{-1.5cm} \sup_{(p,q) \in E\times E}\binom{n}{2}^{-1} \hspace{-0.4 cm} \sum_{1 \leq i_1 < i_2 \leq n} \left \vert h_{p,q}\left(X_{i_1},X_{i_2}\right) - E \left[ h_{p,q}\left(X_{i_1},X_{i_2}\right) \right] \right\vert  \ccs 0 
\end{align*}
when $n \rightarrow \infty$, and from the dominated convergence theorem we derive 
\begin{small}
\begin{align*}
 E \left[ h_{p,q}\left(X_{1},X_{2}\right) \right] \leq  2C \nu \left[ B\left(p,R^*\right)  \triangle B\left(q,R^*\right) \right] \defeq \gamma(\epsilon) \rightarrow 0, \ \mbox{as}  \ \ \epsilon \rightarrow 0.
\end{align*} 
\end{small}
\QEDB
\end{enumerate}
\item \textbf{Proof of Theorem \ref{consistenciamanifold}.}

Given $\epsilon >0$, from Theorem  \ref{desvanecimientoinfinitomanifold}, for $M$ and  $n$ sufficiently large we have that  
\begin{align*}
\sup_{d(x,q) \geq M } \vert \widehat{BD}_n(x) - BD(x) \vert \leq \sup_{d(x,q) \geq M } \vert \widehat{BD}_n(x)\vert + \\
 &\hspace{-7cm} + \sup_{d(x,q) \geq M } \vert BD(x) \vert < \epsilon/3 \quad a.s.
\end{align*}
Given the compact set $\bar{B}(q,M)$,  we can make $\delta>0$ small enough so that we can use Theorem  \ref{teoremacontinuidadmanifold}. Since $\M$ is separable, there exists a finite set
$D\defeq \{x_1, \ldots x_m\} \subset \bar{B}(q,M) $ such that $\cup_{j=1}^m B(x_j, \delta) \supset \bar{B}(q,M)$. For $x \in B(x_j, \delta)$,
$$\vert \widehat{BD}_n(x) - BD(x) \vert \leq \vert \widehat{BD}_n(x) - \widehat{BD}_n(x_i) \vert  + \vert \widehat{BD}_n(x_i) - BD(x_i) \vert $$
$$+ \vert BD(x_i) - BD(x) \vert  \leq 
 2\epsilon + R_n + \sup_{y \in D} \vert \widehat{BD}_n(y) - BD(y) \vert.$$
We can conclude that
$$\sup_{x \in \bar{B}(q,M)} \left \vert \widehat{BD}_n(x) - BD(x) \right \vert \leq 2\epsilon + R_n +\sup_{y \in D} \vert \widehat{BD}_n(y) - BD(y) \vert.$$

On the other hand, since the kernel of the order $2$ $U$-statistic  $\widehat{BD}_n(p)$ is bounded between $0$ and $1$, we can use the Hoeffding inequality   obtaining that for $s>0$, 
\begin{displaymath}
P\left( \vert \widehat{BD}_n(y) - BD(y) \vert> s \right) \leq  2 \exp \left\{ - \frac{ns^2}{2} \right\}, \textrm{for $n>2$}.
\end{displaymath}
Using this inequality, we get 

\begin{align*}
P \Big( \sup_{y \in D} \vert \widehat{BD}_n(y) - BD(y) \vert > s \Big)  \leq  m P \Big(\vert \widehat{BD}_n(y) - BD(y) \vert > s \Big)  \\
&  \hspace{-5 cm}  \leq 2 m \exp \left\{ - \frac{ns^2}{2} \right\},
\end{align*}
for $n>2$. Finally, the Borel--Cantelli Lemma entails that $$\sup_{y \in D} \vert \widehat{BD}_n(y) - BD(y) \vert$$ converges a.s. to 0, which together with the result in the first step concludes the proof. \QEDB 

\item \textbf{Proof of Theorem \ref{normalidad}.}
By assumption, we have that $\mathcal D$ (defined in \ref{estimadorprofundidad_euc})  has a finite $VC$--dimension.  Lastly, this implies that the assumptions in Proposition $10$ in  \cite{gine1996} hold for a order $2$ $U$-statistic and the asymptotic distribution is derived from  Theorem 4.10 in  \cite{arcones1993}.
\QEDB

\item \textbf{Proof of Theorem \ref{hilb}.}

We consider  the Hilbert product space $\mathcal{H}^2 \defeq \mathcal{H \times H}$ and the subclass of functions
$$\mathcal{F}:= \{f_x: x=(x_1,x_2) \in \mathcal{H}^2\} \ \ \mbox{where} \ \ f_x(z) = \I_{ \{ \langle x_1-z, x_2-z \rangle \leq 0 \} }.$$
\\
First, observe that it suffices to prove the theorem for the following V-statistic
\\
 $$\widetilde{BD}_n(z) \defeq \frac{1}{n^2}  \sum_{1 \leq i_1,i_2 \leq n} \I_{B_{ X_{i_1} X_{i_2}}}(z)=  (1 - \frac{1}{n}) \widehat{BD}_n(z).$$  
Obviously,
$$\sup \left \{  \left \vert f_x(z)-f_x(t)  \right \vert : x \in \mathcal{H}^2 \,\, \textrm{and} \,\, z,t  \in \mathcal{H} \right \} \leq 1 < \infty.$$
Next observe  that, for $\epsilon>0$ $$ \left \vert \I_{ \{ \langle t^{'}_1-y, t^{'}_2-y \rangle \leq 0\} }  - \I_{ \{ \langle t_1-y, t_2-y \rangle \leq 0 \} } \right \vert > \epsilon  \Leftrightarrow   \langle t^{'}_1-y, t^{'}_2-y \rangle  \langle t_1-y, t_2-y \rangle \leq 0.$$ 
\\
Let $x \defeq (x_1,x_2) \in \mathcal{H}^2$, $\Vert (t_1, t_2) \Vert_{\mathcal{H}^2} \defeq  \max \left \{    \Vert t_1 \Vert_{\mathcal{H}} ,  \Vert t_2 \Vert_{\mathcal{H}}    \right \}$ and $(t_1,t_2)$, $(t^{'}_1,t^{'}_2)  \in B \left( (x_1,x_2), \delta \right)$  the ball in $\mathcal{H}^2$ centred at $(x_1,x_2)$ and radius $\delta$ with respect to the norm $  \Vert . \Vert_{\mathcal{H}^2}$.  We just consider the case where  $\langle t^{'}_1-y, t^{'}_2-y \rangle \geq 0$ and $\langle t_1-y, t_2-y \rangle \leq 0$ since the case where
  $\langle t^{'}_1-y, t^{'}_2-y \rangle \leq 0$ and $\langle t_1-y, t_2-y \rangle \geq 0$ follows analogously.
 \\
We have that,  $\frac{t_1 + t_2}{2}$ and  $\frac{t^{'}_1+t^{'}_2}{2}$ are in $B \left (\frac{x_1+x_2}{2}, \delta \right)$. In a Hilbert  space   $\mathcal{H}$ for all $y,t_1,t_2 \in \mathcal{H}$  we have that $y \in B_{ t_1 t_2} \Leftrightarrow  \langle t_1 - y, t_2 - y \rangle \leq 0$,  therefore 
 $$\quad y \in B \left(\frac{t_1+t_2}{2}, \frac{ \Vert t_1- t_2 \Vert}{2}\right) \cap  B^c\left(\frac{t^{'}_1+t^{'}_2}{2}, \frac{ \Vert t^{'}_1- t^{'}_2 \Vert}{2}\right).$$ Moreover, we have that,
\\
\begin{itemize}
\item from $\frac{\Vert t_1 -t_2\Vert}{2} \leq  \frac{\Vert t_1 -x_1\Vert}{2} +\frac{\Vert x_1 -x_2\Vert}{2} +\frac{\Vert x_2 -t_2\Vert}{2} \leq                      \frac{\Vert x_1-x_2\Vert}{2} +\delta$ it follows that
$ \Vert y -\frac{x_1+x_2}{2} \Vert \leq \Vert y -\frac{t_1+t_2}{2} \Vert + \Vert \frac{t_1+t_2}{2} -\frac{x_1+x_2}{2} \Vert \leq  \Vert \frac{t_1-t_2}{2}\Vert+ \delta \leq \frac{\Vert x_1-x_2\Vert}{2} +2\delta $. Therefore, 
$$ y \in B \left(\frac{x_1+x_2}{2}, \frac{ \Vert x_1- x_2 \Vert}{2} + 2\delta\right)$$.
\item  analogously, from $\frac{\Vert t^{'}_1 -t^{'}_2\Vert}{2} \geq \frac{\Vert x_1 - x_2\Vert}{2} - \delta$ it follows that $$y \notin B \left(\frac{x_1+x_2}{2}, \frac{ \Vert x_1- x_2 \Vert}{2} -2 \delta\right),$$
for $\delta$ small enough.
\end{itemize} 
Next take $K_\gamma  \in \mathcal{H} $ compact 
et such that, if we denoted $X=(X_1,X_2)$,   $P\left( X \in  (K_\gamma \times K_\gamma)^c  \right)< \gamma$. 
\\
Let $M =  \frac{X_1+X_2}{2}$ and   $R=\frac{\Vert X_1- X_2 \Vert }{2}$,  we have that
 \begin{align*}
  \lim_{\delta \rightarrow 0} \sup_y P \left(  \sup_{t, t^{'} \in B(x, \delta)} \left \vert \I_{ \{ \langle t^{'}_1-y, t^{'}_2-y \rangle \leq 0\} }  - \I_{ \{ \langle t_1-y, t_2-y \rangle \leq 0 \} } \right\vert >\epsilon \right)  \\
  & \hspace{-11cm} \leq   \lim_{\delta \rightarrow 0} \sup_y P \left( y \in   B(M, R+ 2\delta)   \setminus  B(M, R- 2\delta) \right)   \\
& \hspace{-11cm}    \leq   \lim_{\delta \rightarrow 0} \sup_y P \left( y \in   B(M, R+ 2\delta)   \setminus  B(M, R- 2\delta) , X \in  K_{\gamma} \times K_{\gamma} \right) +\gamma. 
  \end{align*}  
 \\
Since $K_\gamma$ is compact, there exists  $L>0$ with   $K_\gamma \subset B(0,L)$, and  let a compact set $K_{\gamma,L}$   such that   $K_{\gamma} \subset K_{\gamma,L} \ominus  B(0,2L+1)$ 
where $K_{\gamma,L} \ominus  B(0,2L+1) \defeq  \left \{ z \in  K_{\gamma,L} : d \left( z,  K^c_{\gamma,L} \right) > 2L +1 \right \}$.
\\
 Let 
$\psi_{\delta}(y) = P \left( y \in   B(M, R+ 2\delta)   \setminus  B(M, R- 2\delta) , X \in  K_{\gamma} \times K_{\gamma} \right)$. 
If $y \in K^c_{\gamma,L}$  then  
$d \left(\frac{x_1+x_2}{2},y \right) \geq d \left(y,x_1 \right) - d\left(x_1,\frac{x_1+x_2}{2}\right) \geq L+1 \geq L + 2\delta \geq  \frac{\Vert x_1 - x_2 \Vert}{2} + 2\delta, $ 
for all $x_1, x_2  \in K_{\gamma}$. which implies that $\psi_{\delta}(y)=0$ for all $y \in K^c_{\gamma,L}$. Therefore we have that
$$\sup_{y \in \mathcal{H}} \psi_{\delta}(y) =  \sup_{y \in K_{\gamma,L}}  \psi_{\delta}(y) =  \psi_{\delta}(y^*_{\delta}), $$ 
for all $\gamma >0$ with $y^*_{\delta}$  being the point where the maximum is attained in the compact set  $K_{\gamma,L}$. 
\\  
The proof will be completed if we show that for all fixed $\gamma$  $$\lim_{\delta \to 0} \psi_{\delta}(y^*_{\delta})=0.
$$ If this is not the case there exists $\eta>0$, $y_n \in  K_{\gamma,L}$ and $\delta_n\rightarrow 0$ such that $$\psi_{\delta_n}(y_n)=P \Big( M \in   B(y_n, R + 2\delta_n)   \setminus  B(y_n, R - 2\delta_n), X \in  K_{\gamma} \times K_{\gamma} \Big ) >\eta.$$ Since $ [0,1] \times K_{\gamma,L}$ is compact there exists a convergence subsequence, which we will denote $(\delta_n,y_n)$,  such that $ (\delta_n,y_n) \rightarrow (0,y)$ for some $ y\in K_{\gamma,L}$. 
Since  $$\psi_{\delta_n}(y_n)=\mathbb E\left[ \mathbb E\left(  I_{\left \{B(y_n, R + 2\delta_n)   \setminus  B(y_n, R- 2\delta_n) \right \}} (M)  I_{ \{  K_{\gamma} \times K_{\gamma} \}} (X) \big \vert R=r \right) \right],$$
then, from assumption \bf HB \rm, the dominated convergence theorem and the continuity of the function  $\psi_{\delta}(y)$ we get that 
\begin{align*}
P \Big( M \in   B(y_n, R + 2\delta_n)   \setminus  B(y_n, R- 2\delta_n), X \in  K_{\gamma} \times K_{\gamma} \Big) \to \\
& \hspace{-10 cm} \to \mathbb E\left[  P \big( M \in B(y, R ), X \in  K_{\gamma} \times K_{\gamma}  \big \vert R=r \big ) \right]=0,
\end{align*}
    which leads to a contradiction.

\QEDB
\end{enumerate}


\begin{thebibliography}{99}

\bibitem[Anderson, 1955]{anderson1955}
Anderson, T.~W. (1955).
\newblock The integral of a symmetric unimodal function over a symmetric convex
  set and some probability inequalities.
\newblock {\em Proceedings of the American Mathematical Society},
  6(2):170--176.

\bibitem[Arcones and Gine, 1993]{arcones1993}
Arcones, M.~A. and Gine, E. (1993).
\newblock Limit theorems for U--processes.
\newblock {\em The Annals of Probability}, 21(3):1494--1542.

\bibitem[Arnaudon et~al., 2013]{barbaresco2013}
Arnaudon, M., Barbaresco, F., and Yang, L. (2013).
\newblock Riemannian medians and means with applications to radar signal
  processing.
\newblock {\em IEEE Journal of Selected Topics in Signal Processing},
  7(4):595--604.

\bibitem[Barnett, 1976]{barnett1976}
Barnett, V. (1976).
\newblock The ordering of multivariate data.
\newblock {\em Journal of the Royal Statistical Society. Series A (General)},
  139(3):318--355.

\bibitem[Bhattacharya and Bhattacharya, 2008]{bhattacharya2008}
Bhattacharya, A. and Bhattacharya, R. (2008).
\newblock {\em Nonparametric statistics on manifolds with applications to shape
  spaces}, 3:282--301.
\newblock Institute of Mathematical Statistics, Beachwood, Ohio, USA.



\bibitem[Billingsley, Patrick and Tops{\o}e, Flemming, 1967]{billingsley1967}
Billingsley, Patrick and Tops{\o}e, Flemming. (1967).
\newblock {\em Uniformity in weak convergence}, volume~7, pages 1--16.
\newblock {Zeitschrift f{\"u}r Wahrscheinlichkeitstheorie und verwandte Gebiete.}

\bibitem{carrizosa1996}
Carrizosa, E. (1996).
\newblock A characterization of halfspace depth.
\newblock {\em Journal of multivariate analysis}, 58(1):21--26.

\bibitem[Chen et~al., 2017]{chen2017}
Chen, M., Gao, C., and Ren, Z. (2017).
\newblock Robust covariance and scatter matrix estimation under huber's
  contamination model.
\newblock {\em To appear in Annals of Statistics}.

\bibitem[Croux et~al., 2017]{ggordaliza2017}
Croux, C., Garcia-Escudero, L.A., Gordaliza, A., Ruwet, C., and Martin, R.
  (2017).
\newblock Robust principal component analysis based on trimming around affine
  subspaces.
\newblock {\em Statistica Sinica}, 27(3):1437--1459.

\bibitem{cuevas2014}
Cuevas, A. (2014).
\newblock A partial overview of the theory of statistics with functional data.
\newblock {\em Journal of Statistical Planning and Inference}, 147:1--23. 

\bibitem[Cuevas and Fraiman, 2009]{fraiman2009}
Cuevas, A. and Fraiman, R. (2009).
\newblock {On depth measures and dual statistics. A methodology for dealing
  with general data.}
\newblock {\em J. Multivariate Analysis}, 100(4):753--766.

\bibitem[Do~Carmo and Flaherty~Francis, 1992]{docarmo1992}
Do~Carmo, M.~P. and Flaherty~Francis, J. (1992).
\newblock {\em Riemannian geometry}, volume 115.
\newblock Birkh{\"a}user Boston.

\bibitem[Donoho and Gasko, 1992]{donoho1992}
Donoho, D.~L. and Gasko, M. (1992).
\newblock Breakdown properties of location estimates based on halfspace depth
  and projected outlyingness.
\newblock {\em The Annals of Statistics}, pages 1803--1827.

\bibitem[Dudley, 1978]{dudley1978}
Dudley, R.~M. (1978).
\newblock Central limit theorems for empirical measures.
\newblock {\em The Annals of Probability}, pages 899--929.

\bibitem[D{\"u}mbgen, 1992]{dumbgen1992}
D{\"u}mbgen, L. (1992).
\newblock Limit theorems for the simplicial depth.
\newblock {\em Statistics \& Probability Letters}, 14(2):119--128.


\bibitem{elmore2006}
Elmore, R.T., Hettmansperger, T.P., and  Xuan, F. (2006)
\newblock Spherical data depth and a multivariate median.
\newblock {\em DIMACS Series in Discrete Mathematics and Theoretical Computer
  Science}, 72:87.

\bibitem{ferri2008}
Ferri, M. and Frosini, P. (2008).
\newblock VC-dimension on manifolds: a first approach.
\newblock {\em Mathematical methods in the applied sciences}, 31(5):589--605.


\bibitem[Fletcher et~al., 2009]{fletcher2009}
Fletcher, P.~T., Venkatasubramanian, S., and Joshi, S. (2009).
\newblock The geometric median on riemannian manifolds with application to
  robust atlas estimation.
\newblock {\em NeuroImage}, 45(1):S143--S152.

\bibitem[Folland, 2013]{folland2013}
Folland, G.~B. (2013).
\newblock {\em Real analysis: modern techniques and their applications}.
\newblock John Wiley \& Sons.

\bibitem[Gin{\'e}, 1996]{gine1996}
Gin{\'e}, E. (1996).
\newblock Empirical processes and applications: an overview.
\newblock {\em Bernoulli}, 2(1):1--28.



\bibitem{hampel1971}
Frank~R Hampel.(1971)
\newblock A general qualitative definition of robustness.
\newblock {\em The Annals of Mathematical Statistics}, 42(6):1887--1896.

\bibitem[Hoffmann-J{\o}rgensen, 1991]{hoffmann1991}
Hoffmann-J{\o}rgensen, J. (1991).
\newblock {\em Stochastic processes on Polish spaces}.
\newblock Number~39. Aarhus Universitet. Matematisk Institut.

\bibitem[Lin et~al., 2016]{lin2016}
Lin, L., Thomas, B.~S., Zhu, H., and Dunson, D.~B. (2016).
\newblock Extrinsic local regression on manifold-valued data.
\newblock {\em Journal of the American Statistical Association to appear}.

\bibitem[Liu, 1990]{liu1990}
Liu, R.~Y. (1990).
\newblock {On a Notion of Data Depth Based on Random Simplices}.
\newblock {\em The Annals of Statistics}, 18(1):405--414.

\bibitem[Liu, 1992]{liu1992}
Liu, R.~Y. (1992).
\newblock Data depth and multivariate rank tests.
\newblock {\em L1-statistical analysis and related methods (Y. Dodge, ed.)},
  pages 279--294.

\bibitem[Liu et~al., 1999]{liu1999}
Liu, R.~Y., Parelius, J.~M., and Singh, K. (1999).
\newblock Multivariate analysis by data depth: Descriptive statistics, graphics
  and inference.
\newblock {\em The Annals of Statistics}, 27(3):783--840.

\bibitem[Liu and Singh, 1992]{liu1992b}
Liu, R.~Y. and Singh, K. (1992).
\newblock Ordering directional data: Concepts of data depth on circles and
  spheres.
\newblock {\em Ann. Statist.}, 20(3):1468--1484.

\bibitem{liu2011}
Zhenyu Liu and Reza Modarres.
\newblock Lens data depth and median.
\newblock {\em Journal of Nonparametric Statistics}, 23(4):1063--1074, 2011.


\bibitem[Mardia et~al., 2008]{mardia2008}
Mardia, K.~V., Hughes, G., Taylor, C.~C., and Singh, H. (2008).
\newblock A multivariate von mises distribution with applications to
  bioinformatics.
\newblock {\em Canadian Journal of Statistics}, 36(1):99--109.

\bibitem[Mardia and Jupp, 2000]{mardia2000}
Mardia, K.~V. and Jupp, P.~E. (2000).
\newblock {\em Directional statistics}.
\newblock John Wiley \& Sons.

\bibitem[Mardia and Voss, 2014]{mardia2014}
Mardia, K.~V. and Voss, J. (2014).
\newblock Some fundamental properties of a multivariate von mises distribution.
\newblock {\em Communications in Statistics-Theory and Methods},
  43(6):1132--1144.

\bibitem[Maronna et~al., 2006]{maronna2006}
Maronna, R., Martin, R.~D., and Yohai, V. (2006).
\newblock {\em Robust statistics}.
\newblock John Wiley \& Sons, Chichester. ISBN.

\bibitem[Moakher, 2005]{moakher2005}
Moakher, M. (2005).
\newblock A differential geometric approach to the geometric mean of symmetric
  positive-definite matrices.
\newblock {\em SIAM Journal on Matrix Analysis and Applications},
  26(3):735--747.

\bibitem[Oja, 1983]{oja1983}
Oja, H. (1983).
\newblock Descriptive statistics for multivariate distributions.
\newblock {\em Statistics \& Probability Letters}, 1(6):327--332.

\bibitem[Patrangenaru and Ellingson, 2015]{patrangenaru2015}
Patrangenaru, V. and Ellingson, L. (2015).
\newblock {\em Nonparametric statistics on manifolds and their applications to
  object data analysis}.
\newblock CRC Press.

\bibitem{pelletier2006}
Pelletier, B. (2006).
\newblock Non-parametric regression estimation on closed riemannian manifolds.
\newblock {\em Journal of Nonparametric Statistics}, 18(1):57--67. 

\bibitem[Pennec, 2006]{pennec2006}
Pennec, X. (2006).
\newblock Intrinsic statistics on riemannian manifolds: Basic tools for
  geometric measurements.
\newblock {\em Journal of Mathematical Imaging and Vision}, 25(1):127--154.

\bibitem[Petersen, 2006]{petersen2006}
Petersen, P. (2006).
\newblock {\em Riemannian geometry}, volume 171.
\newblock Springer.


\bibitem{pizer2017}
Pizer, S.M.  and marron J.S. (2017).
\newblock Object statistics on curved manifolds.
\newblock In {\em Statistical Shape and Deformation Analysis}, pages 137--164.
  Elsevier.




\bibitem[Pollard, 2012]{pollard2012}
Pollard, D. (2012).
\newblock {\em Convergence of stochastic processes}.
\newblock Springer Science \& Business Media.





\bibitem[Serfling, 1980]{serfling1980}
Serfling, R. (1980).
\newblock {\em {Approximation theorems of mathematical statistics}}.
\newblock {Wiley series in probability and mathematical statistics :
  Probability and mathematical statistics}. Wiley, New York edition.


\bibitem{serfling2006}
Serfling, R. (2006).
\newblock Depth functions in nonparametric multivariate inference.
\newblock {\em DIMACS Series in Discrete Mathematics and Theoretical Computer
  Science}, 72:1.
  
  
\bibitem[Serfling and Zuo, 2000]{serfling2000}
Serfling, R. and Zuo, Y. (2000).
\newblock {General notions of statistical depth function}.
\newblock {\em The Annals of Statistics}, 28(2):461--482.


\bibitem{shahsavarifar2018}
Shahsavarifar, R. and Bremner, D. (2018).
\newblock Computing the planar $beta$-skeleton depth.
\newblock {\em preprint arXiv:1803.05970.} 


\bibitem[Steele, 1975]{steele1975}
Steele, J.~M. (1975).
\newblock {\em Combinatorial entropy and uniform limit laws}.
\newblock Department of Mathematics, Stanford University.

\bibitem[Szabados, 1989]{szabados1989}
Szabados, T. (1989).
\newblock {On the Glivenko-Cantelli theorem for balls in metric spaces}.
\newblock {\em Studia Scientiarum Mathematicarum Hungarica}, 24:473--481.


  
\bibitem[Tukey, 1975]{tukey1975}
Tukey, J.~W. (1975).
\newblock Mathematics and the picturing of data.
\newblock In {\em Proceedings of the international congress of mathematicians},
  volume~2, pages 523--531.

\bibitem[Vardi and Zhang, 2000]{vardi2000}
Vardi, Y. and Zhang, C.-H. (2000).
\newblock The multivariate l1-median and associated data depth.
\newblock {\em Proceedings of the National Academy of Sciences},
  97(4):1423--1426.

\bibitem[Zhu et~al., 2009]{zhu2009}
Zhu, H., Chen, Y., Ibrahim, J.~G., Li, Y., Hall, C., and Lin, W. (2009).
\newblock Intrinsic regression models for positive-definite matrices with
  applications to diffusion tensor imaging.
\newblock {\em Journal of the American Statistical Association},
  104(487):1203--1212.
\end{thebibliography}
\end{document}